\documentclass{article}

\usepackage{amsmath,amssymb,amsthm,latexsym,hyperref,url,graphicx,cite}

\newtheorem{theorem}{Theorem}

\newcommand{\RR}[2]{R_{#1}\left(#2\right)}

\newcommand{\mytexttilde}{\raise.17ex\hbox{$\scriptstyle\sim$}}

\newcommand{\Z}{\mathbb{Z}}

\newcommand{\ra}{\rightarrow}
\newcommand{\E}{\mathcal{E}}

\begin{document}

\title{A Note on a Question of Erd\H{o}s \& Graham}
\author{Kellen Myers}
\date{\today}

\maketitle
\begin{abstract}
Erd\H{o}s \& Graham ask whether the equation $x^2+y^2=z^2$ is partition regular,
i.e. whether it has a finite Rado number. This note provides a lower bound and also
states results in the affirmative for two similar quadratic equations.
\end{abstract}

\section{Introduction}

In \cite{Schur}, Schur shows that for any finite coloring of the positive integers (a function $\chi:\Z^+\ra\{1,2,\dots,r\}$), there will be a triple $(x,y,z)$ such that $x+y=z$ and $\chi(x)=\chi(y)=\chi(z)$. This triple is said to be a monochromatic solution to the equation $x+y=z$. The least $N$ such that this statement holds for any coloring $\chi:\{1,2,\dots,N\}\ra\{1,2,\dots,r\}$ is called the $r$-color Schur number and is denoted $S(r)$.

It is easy to see that $S(2)=5$, and one might note trivially that $S(1)=2$. The values of $S(3)$ and $S(4)$ are also known. We can associate such a quantity to any equation $\E$, not just $x+y=z$, which denote $\RR{r}{\E}$. In cases where no such $N$ exists, we say $\RR{r}{\E}=\infty$. An equation $\E$ is called $r$-regular if the quantity $\RR{r}{\E}$ is finite, and it is called regular if it is $r$-regular for all $r$.

In \cite{Rado}, Rado provides necessary and sufficient conditions for linear, homogeneous equations to be regular. He also gives necessary and sufficient conditions (essentially non-triviality) for such equations to be 2-regular. For that reason, we call these quantities ``Rado numbers." There have been many papers, starting with \cite{BB82}, giving certain Rado numbers (often parametrized families of equations). In most cases the equations are linear and there are 2 colors.

In this note, we state three results from a forthcoming paper \cite{MP14} that do not fall into these categories. The equations are quadratic, and in one case, the number of colors is 3.

In \cite{EG80}, Erd\H{o}s and Graham ask whether the equation $x^2+y^2=z^2$ is 2-regular. In \cite{Gra08}, Graham notes that it is not clear which answer is correct. In \cite{FH13} it is proved that $9x^2+16y^2=n^2$ (along with a family of similar quadratic equations) is 2-regular, but only $x$ and $y$ are supposed to be monochromatic (note $n$). We offer the following three results, which have been a part of ongoing work to settle the question:

{\theorem $\RR{2}{x^2+y^2+z^2=w^2}=105$.}

{\theorem $\RR{3}{x_1^2+x_2^2+x_3^2+x_4^2=y_1^2+y_2^2+y_3^2}=32$.}

{\theorem $\RR{2}{x^2+y^2=z^2}>6500$.}

In the third of these results, we do not exclude the possibility that it is infinite.

These results are all obtained computationally and will be detailed in \cite{MP14}. The first result also inspires a new sequence representing the Rado number for the equation $x_1^2 + x_2^2 + \dots + x_k^2 = z^2$, which is tabulated as follows:

$$\begin{array}{|c|c|c|c|c|c|c|c|c|c|c|c|c|c|c|c|c|}\hline
k = & 2 & 3   & 4  & 5  & 6  & 7  & 8  & 9  & 10 & 11 & 12 & 13 & 14 & 15 & 16 & 17\\ \hline
N=  & ? & 105 & 37 & 23 & 18 & 20 & 20 & 15 & 16 & 20 & 23 & 17 & 21 & 26 & 17 & 23\\ \hline
\end{array}$$

This is now entry \href{http://oeis.org/A250026}{A250026} in the Online Encyclopedia of Integer Sequences.

\bibliographystyle{amsalpha}
\bibliography{SRN}

\end{document}